\documentclass[11pt]{article}

\usepackage{amsmath}
\usepackage{amsthm}
\usepackage{amsfonts}
\usepackage{amssymb}
\usepackage{mathtools}

\usepackage{hyperref}
\usepackage[compress,sort]{cite}

\usepackage{cases}
\usepackage{xcolor}

\newcommand{\dz}{\partial_z}
\newcommand{\dt}{\partial_t}
\newcommand{\ddh}{\partial_h}
\newcommand{\dvh}{\mathrm{div}_h\,}
\newcommand{\crlh}{\mathrm{curl}_h\,}
\newcommand{\nablah}{\nabla_h}
\newcommand{\Deltah}{\Delta_h}
\newcommand{\idx}{\, d \vec x}

\newcommand{\norm}[2]{\bigl\Arrowvert #1 \bigr\Arrowvert_{#2}}
\newcommand{\Lnorm}[1]{L^{#1}(\Omega)}
\newcommand{\Hnorm}[1]{H^{#1}(\Omega)}

\newcommand{\mrm}[1]{\mathrm{#1}}

\theoremstyle{plain}
\newtheorem{thm}{Theorem}[section]
\newtheorem{lm}{Lemma}

\theoremstyle{definition}

\theoremstyle{remark}
\newtheorem{remark}{Remark}

\numberwithin{equation}{section}
\allowdisplaybreaks
\title{Derivation of a generalized quasi-geostrophic approximation for inviscid
flows in a channel domain: The fast waves correction
}

\author{Claude Bardos
	\footnote{Laboratoire J.-L. Lions, BP187, 75252 Paris Cedex 05, France. \href{mailto:claude.bardos@gmail.com}{claude.bardos@gmail.com}
	}, 
\quad 
Xin Liu
	\footnote{Department of Mathematics, Texas A\&M University, College Station, TX 77843-3368, USA
	\href{mailto:stleonliu@gmail.com}{stleonliu@gmail.com}
	}, 
\quad and \quad 
Edriss S. Titi
	\footnote{Department of Mathematics, Texas A\&M University, College Station, TX 77843-3368, USA; Department of Applied Mathematics and Theoretical Physics, University of Cambridge, Cambridge CB3 0WA UK; also Department of Computer Science and Applied Mathematics, Weizmann Institute of Science, Rehovot 76100, Israel.
	\href{mailto:titi@math.tamu.edu}{titi@math.tamu.edu} and \href{mailto:Edriss.Titi@damtp.cam.ac.uk}{Edriss.Titi@damtp.cam.ac.uk}
	}
}

\begin{document}
\maketitle

\begin{abstract}
    This paper is devoted to investigating the rotating Boussinesq equations of inviscid, incompressible flows with both fast Rossby waves and fast internal gravity waves. The main objective is to establish a rigorous derivation and justification of a new generalized quasi-geostrophic approximation in 
    a channel domain with no normal flow at the upper and lower solid boundaries, taking into account the resonance terms due to the fast and slow waves interactions. Under these circumstances, We are able to obtain uniform estimates and compactness without the requirement of either well-prepared initial data (as in \cite{bourgeoisValidityQuasigeostrophicModel1994}) or domain with no boundary (as in \cite{embidAveragingFastGravity1996}). 
    In particular, the nonlinear resonances and the new limit system, which takes into account the fast waves correction to the slow waves dynamics, are also identified without introducing Fourier series expansion.
    The key ingredient includes the introduction of (full) generalized potential vorticity.
    
    \bigskip
    
    {\noindent\bf Keyworks: } Quasi-Geostrophic approximation, singular limit, Rossby waves, internal gravity waves, bounded domain, fast-slow waves interaction, potential vorticity.
    
    \bigskip
    
    
\end{abstract}

\tableofcontents

\section{Introduction}\label{sec:introduction}

We consider an inviscid incompressible fluid in a periodic channel domain $ \Omega := \Omega_h \times (0,h) \subset \mathbb R^3 $, with horizontal periodic domain $ \Omega_h := \mathbb T^2 = (0,1)^2 $ and vertical domain height $ h \in (0,\infty) $. Denote by $ v \in \mathbb R^2 $ the horizontal velocity, $ w \in \mathbb R $ the vertical velocity, $ p \in \mathbb R $ the pressure, and $ \rho \in \mathbb R $ the density, respectively. Let the following be the typical characteristic physical scales for length, time, velocity, density, and pressure:
\begin{equation*}
    \begin{aligned}
        & L && \text{length scale} \\
        & U && \text{mean advective velocity} \\
        & T_e := \dfrac{L}{U} && \text{eddy trunover time} \\
        & T_R := f^{-1} && \text{rotation time} \\
        & \rho_b && \text{mean density} \\
        & \overline p && \text{mean pressure.}
    \end{aligned}
\end{equation*}
Furthermore, set $ \overline \rho = \overline\rho(z) $ to be the background density stratification, which is assumed to be linear in the vertical coordinate, and decompose the density into the sum of stratification $ \overline \rho $ and deviation $ \rho_b\theta $, i.e.,  
\begin{equation*}
    \rho = \rho_b \theta + \overline \rho
\end{equation*}
The buoyancy (Brunt-V\"ais\"al\"a) frequency is defined as
\begin{equation*}
    N:= \biggl(- \dfrac{g \partial_z \overline \rho}{\rho_b} \biggr)^{1/2},
\end{equation*}
and the corresponding buoyancy time scale is
\begin{equation*}
    T_N := N^{-1}.
\end{equation*}
In this geophysical situation, one can introduce the following relevant non-dimensional numbers:
\begin{align*}
    &\text{the Rossby number} && \mathrm{Ro}:= 
    \dfrac{U}{Lf} \\
    &\text{the Froude number} &&  \mathrm{Fr}:= 
    \dfrac{U}{LN} \\
    &\text{the Euler number} && \overline{P}:= \dfrac{\overline p}{\rho_bU^2} \\
    & &&  \Gamma := \dfrac{gL}{U^2},
\end{align*}
see, e.g., \cite{MajdaAtmosphereOcean}.
With such notations, the dimensionless rotating Boussinesq equations are given by
\begin{subequations}\label{sys:rt-bn-eq}
    \begin{align}
    \label{eq:bn_ns_v}
        \partial_t v +  v \cdot \nablah v + w \dz v  + \dfrac{1}{\mathrm{Ro}} v^\perp  + \overline{P} \nablah p  = 0, \\
    \label{eq:bn_ns_w}
        \dt w + v \cdot \nablah w + w \dz w  + \overline{P}\dz p  - \Gamma\theta = 0, \\
    \label{eq:bn_ns_tm}
        \dt \theta + v \cdot \nablah \theta + w \dz \theta  + \dfrac{1}{\Gamma \cdot \mathrm{Fr}^2} w  = 0, \\
    \label{eq:bn_ns-dv}
        \dvh v + \dz w  = 0, 
\end{align}
with 
\begin{align}
	w \vert_{z=0,h} = 0 & \qquad \text{i.e., the impermeable boundary condition,} \label{bc:bn_ns}
\end{align}
\end{subequations}
see, e.g., \cite{MajdaAtmosphereOcean}.

\bigskip 

In this paper, we consider the quasi-geostrophic scale where
\begin{itemize}
    \item The Rossby number is small $$ \mathrm{Ro} = \varepsilon \ll 1; $$
    \item The flow is in geostropic balance, i.e., the rotation and the pressure forces are in balance,
    $$ \overline P = \dfrac{1}{\mathrm{Ro}}; $$
    \item The Froude number is small and equal to the Rossby number, $$ \mathrm{Fr} = \mathrm{Ro}; $$
    \item The non-dimensional number $ \Gamma $ is in balance with the inverse of the Froude number $$ \Gamma = \dfrac{1}{\mathrm{Fr}}. $$
\end{itemize}
Then the rotating Boussinesq equations \eqref{sys:rt-bn-eq} become
\begin{subequations}\label{sys:ns_rt_by}
\begin{align}
    \label{eq:ns_v}
        \partial_t v +  v \cdot \nablah v + w \dz v + \dfrac{1}{\varepsilon} v^\perp + \dfrac{\nablah p}{\varepsilon} = 0, \\
    \label{eq:ns_w}
        \dt w + v \cdot \nablah w + w \dz w  + \dfrac{\dz p}{\varepsilon}  - \dfrac{\theta}{\varepsilon}  = 0, \\
    \label{eq:ns_tm}
        \dt \theta + v \cdot \nablah \theta + w \dz \theta  + \dfrac{w}{\varepsilon}  = 0, \\
    \label{eq:ns-dv}
        \dvh v + \dz w = 0, 
\end{align}
with 
\begin{align}
	w \vert_{z=0,h} = 0. 
 \label{bc:ns}
\end{align}
\end{subequations}
We refer the reader to \cite[section 7.4]{MajdaAtmosphereOcean} for the detailed derivation of system \eqref{sys:ns_rt_by}. We remark that, the small Rossby number, i.e. $\mathrm{Ro}\ll1 $, induces the fast Rossby waves, and the small Froude number, i.e. $ \mathrm{Fr} \ll 1 $, induces the fast internal gravity waves. In our setting, i.e., system \eqref{sys:ns_rt_by}, both Rossby and gravity waves are fast and they are coupled. In particular, they have the same scale. 

\bigskip

The goal of this work is to investigate the asymptotic limit of system \eqref{sys:ns_rt_by} as $ \varepsilon \rightarrow 0^+ $ in the channel domain $ \Omega $., i.e., the quasi-geostrophic approximation, taking into account the fast-slow waves interaction and their corresponding resonance terms. 

Similar problem has been studied in the case of ``well-prepared'' initial data by Bourgeois and Beale in \cite{bourgeoisValidityQuasigeostrophicModel1994}, where the convergence, as well as the convergence rate, of solutions to that of quasi-geostrophic equations (\eqref{eq:ptl-vt} and \eqref{eq:sl-bc}, below) is proved. In particular, the well-prepared initial data are chosen so that there are only slow waves in the dynamics and no contribution of the fast waves. That is, the initial data is close to the geostrophic balance (see \eqref{eq:gtph-v}--\eqref{eq:gtph-tm}, below). We remark that \cite{bourgeoisValidityQuasigeostrophicModel1994} assumes that $ \dz p^0\vert_{z=0,h} = 0 $ together with the balanced initial data. 
This guarantees that the system of equations satisfy some symmetry, and eventually can be extended periodically to a system into $ \mathbb T^3 $, i.e., there is no boundary effect as if one has a virtual boundary. The general convergence theory when $ \dz p^0\vert_{z=0,h} \neq 0 $ is still open. Here $ p^0 $ is the stream function associated with the potential vorticity as in \eqref{def:potential-vorticity}. The existence of weak solutions for these quasi-geostrophic equations is established in \cite{novackInviscidThreeDimensional2020, puelGlobalWeakSolutions2015}

Taking into account the fast waves, but without physical boundary (i.e., in $ \mathbb T^3 $), Embid and Majda studied the nonlinear resonances and established the asymptotic limit of system \eqref{sys:ns_rt_by} in \cite{embidAveragingFastGravity1996, embidLowFroudeNumber1998, majdaAveragingFastGravity1998}. The limiting system is the quasi-geostrophic equation \eqref{eq:ptl-vt} with nonlinear resonances on the right-hand side, while the velocity and the temperature in the limiting quasi-geostrophic equations are given by \eqref{eq:gtph-v} and \eqref{eq:gtph-w}, below, respectively. 

In the case with vanishing viscosity, an Ekman boundary layer will arise in the channel domain, which leads to Ekman pumping. This is verified in \cite{desjardinsDERIVATIONQUASIGEOSTROPHICPOTENTIAL1998}, in the case with well-prepared initial data (i.e., slow waves only). To the best of the authors' knowledge, the asymptotic limit taking into account both the fast waves and the Ekman pumping is open. 
The global well-posedness of solutions to the quasi-geostrophic system with Ekman pumping was established in \cite{novackGlobalTimeClassical2018}.

\bigskip

In this paper, we introduce the notion of (full) {\bf generalized potential vorticity} (i.e., $ \Phi $ and $ \Psi $ defined in \eqref{def:pv} and \eqref{def:g_pv}, below, respectively), which allows us to separately describe the slow and the fast waves of the dynamics of system \eqref{sys:ns_rt_by} in a channel domain without introducing any boundary layer. Moreover, the interaction between the slow and fast waves can be easily tracked and investigated. Therefore, we are able to establish the asymptotic limit as $ \varepsilon \rightarrow 0^+ $ in the channel for general initial data. In particular, we drop the requirement of well-prepared initial data or periodic spatial domain required in \cite{bourgeoisValidityQuasigeostrophicModel1994} and \cite{embidAveragingFastGravity1996}, respectively. In addition, the fast waves correction to the slow dynamics is identified as a new resonance term.

We remark that in our context, the terms slow (fast) waves and slow (fast) dynamics, as well as well-prepared (ill-prepared or general) initial data and balanced (unbalanced) initial data are interchangeable, respectively. This terminology is widely used in the literature.

Before stating the main results in detail, we would like to put this work in the context of the study of asymptotic limit in the following subsection. 

\subsection{Asymptotic limit and boundary layer}\label{sec:introduction-asymptotic_limit}

We should stress that the following references are by no mean exhaustive. 

\bigskip

The study of low Mach number limit of the compressible flows was pioneered by Klainerman and Majda in \cite{Klainerman1981, Klainerman1982}, where the convergence with only slow waves (i.e., well-prepared initial data) was shown in domains without boundary. In $ \mathbb R^3 $, Ukai in \cite{Ukai1986} showed the dispersion of the fast acoustic waves and thus established the low Mach number limit with large acoustic waves. As pointed out in \cite{DesjardinsGrenier1999}, such dispersion in $ \mathbb{R}^3 $ is characterized by the Strichartz estimate \cite{Keel1998, strichartzRestrictionsFourierTransforms1977}. 
In the case of $ \mathbb T^3 $, \cite{Lions1998a} showed the weak convergence of low Mach number limit for compressible flows by investigating the nonlinear resonances of fast acoustic waves. The general theory of fast singular limit was developed by Schochet in \cite{Schochet1988, Schochet1994} for hyperbolic systems, which was later extended to parabolic systems in \cite{Gallagher1998}.
We refer the reader to 
\cite{Metivier2001, Alazard2005, Alazard2006, Danchin2002per, Danchin2005, Feireisl2016HandbookSec, Feireisl2007d, Masmoudi2001} and the references therein for more studies of low Mach number limit in domains without boundary. 
When there is physical boundary in the underlying domain, the low Much number limit of viscous flows may give rise to a boundary layer. This is first studied in terms of eigenvalue-eigenfunction pairs in \cite{jiangConstructionBoundaryLayers2015}. Recently in \cite{masmoudiUniformRegularityCompressible2022}, by introducing uniform estimates in the co-normal Sobolev norm, together with some $ L^\infty $ estimates, the low Mach number limit of compressible viscous flows is established in smooth domain with Navier-slip boundary condition and general initial data. However, the corresponding low Mach number limit with no-slip boundary condition is still open. 

\bigskip

Meanwhile, in the vanishing viscosity limit of the incompressible Navier--Stokes equations with no-slip boundary condition, the Prandtl boundary layer was introduced by Prandtl in 1904 \cite{prandtlMotionFluidsVery1928} and became the paradigm of further mathematical studies. See, e.g., \cite{weinanBoundaryLayerTheory2000} for a derivation of the Prandtl equations. However it turned out to be the most singular. The boundary layer is due to the no-slip boundary condition for the Navier-Stokes and since this effect is not present at the level of the Euler equation, a discontinuity appears in the zero viscosity limit. Due to the nonlinearity of the problem such singularity may escape from the boundary layer and propagate in the fluid. This is one of the main source of turbulence, and as a consequence the Prandtl boundary layer is strongly unstable, and therefore may exist only for short time and under strict regularity hypothesis, see, e.g., \cite{Maekawa2014, sammartinoZeroViscosityLimit1998a, sammartinoZeroViscosityLimit1998}.
A direct proof of such asymptotic limit, with the incompressible Euler equations as the limiting equations, without introducing the boundary layer correction can be found in \cite{bardosInviscidLimit2d2021, nguyenInviscidLimitNavier2018}. 
For general, smooth, but not analytic, initial data, the vanishing viscosity limit is still an open challenging problem. The pioneer work in this direction is by Kato \cite{Kato_1984}. See, also,  \cite{Bardos_Titi_2007,Bardos_Titi_2013} and references therein for related results. 


\bigskip

With fast rotation and vanishing viscosity (but no fast internal waves) in a domain with no-slip boundary condition, the Ekman boundary layer may arise, which is an important phenomenon in the atmospheric and oceanic study (see, for instance, \cite{pedloskyGeophysicalFluidDynamics1987, MajdaAtmosphereOcean}). In \cite{grenierEkmanLayersRotating1997} and \cite{masmoudiEkmanLayersRotating2000}, the asymptotic limit of fast rotation and vanishing viscosity with the Ekman boundary layer correction was established for flows with and without fast waves, respectively. 

\bigskip

With only fast rotation in a domain without boundary ($ \mathbb T^3 $ or $ \mathbb R^3 $), the asymptotic limit of the Euler or Navier--Stokes equations was studied in \cite{babinGlobalSplittingIntegrability1996, babinFastSingularOscillating2000, babinGlobalRegularity3D1999}, where the limit dynamics is characterized by two dimensions three components (2D3C) flows, and the prolonging effect of fast rotation on the life-span of the solution was established. Such a regularizing effect of fast rotation was demonstrated in the case of a simple convection model in \cite{Babin_Ilyin_Titi2011,liuRotationPreventsFinitetime2004}.  See also \cite{ghoulEffectRotationLifeSpan2022, linEffectFastRotation2022} for the study in the primitive equations, and \cite{cheminMathematicalGeophysicsIntroduction2006} for some examples in the study of mathematical geophysics, including the aforementioned Ekman boundary layer.

\bigskip

As mentioned before, 
in this paper, we study the singular limit $ \varepsilon \rightarrow 0^+ $ of system \eqref{sys:ns_rt_by} in the periodic channel domain $ \Omega = \mathbb T^2 \times (0,h) $. In particular, it will be established that the fast rotation induced by strong Coriolis force in \eqref{eq:ns_v} suppresses the possible emergence of a boundary layer near the boundary.

\subsection{Main results}\label{sec:introduction-main_results}

The first main result of this paper is the following:

\begin{thm}[Uniform-in-$\varepsilon$ estimate]\label{thm:uniform}
    Consider the initial data $$ (v_\mrm{in}, w_\mrm{in},\theta_\mrm{in}) \in  H^3(\Omega) $$ of the solution $ (v,w,\theta)$ to system \eqref{sys:ns_rt_by}, satisfying the compatibility conditions $ \dvh v_\mrm{in} + \dz w_\mrm{in} = 0 $ and $ w_\mrm{in}\vert_{z=0,h}=0 $. Then there exists $ T, C_\mrm{in} \in(0,\infty) $, depending only on the initial data and independent of $ \varepsilon $, such that
    \begin{equation}\label{thm-est:uniform}
        \sup_{0\leq t \leq T} \norm{v(t),w(t),\theta(t)}{\Hnorm{3}} \leq C_\mrm{in}.
    \end{equation}
\end{thm}
\begin{proof}
The proof of this theorem is done in section \ref{sec:uniform_est}.
\end{proof}

The local well-posedness theory of solutions in $ H^3(\Omega) $ to system \eqref{sys:ns_rt_by} for fixed $ \varepsilon \in (0,1) $ is classical and thus is omitted here. See, for instance, \cite{Kato_Lai_1984}. With continuity arguments, the uniform estimate \eqref{thm-est:uniform} implies the uniform-in-$\varepsilon$ local well-posedness with initial data as in the theorem. 

\bigskip

Our second main result of this paper is to investigate the limit system, as follows:

\begin{thm}[Convergence theory]\label{thm:convergence}
Let $ T > 0 $ be as in Theorem \ref{thm:uniform}, and 
let $ (\Phi, \Psi, H_0, H_h, Z) $ be defined as in \eqref{def:pv}--\eqref{def:bc_v}, below. Then there exists a subsequence of $ \varepsilon $ that as $ \varepsilon \rightarrow 0^+ $, one has the following convergence in strong topology:
\begin{align}
    \Phi & \rightarrow \Phi_p && \text{in} && C([0,T];H^1(\Omega)), \\
    H_0, H_h & \rightarrow H_{p,0}, H_{p,h} && \text{in} && C([0,T];H^{3/2}(\mathbb T^2)), \\
    e^{\mp i \frac{t}{\varepsilon}}(\Psi \pm i \Psi^\perp) & \rightarrow \psi_{p,\pm} && \text{in} && C([0,T];H^1(\Omega)), \\
    \shortintertext{and}
    e^{\mp i \frac{t}{\varepsilon}}(Z \pm i Z^\perp) & \rightarrow z_{p,\pm} && \text{in} && C([0,T];H^2(\Omega)),
\end{align}
and in suitable weak-$*$ topology (see section \ref{sec:convergence_1}), 
the limit 
\begin{equation}
    (\Phi_p, H_{p,0}, H_{p,h}, \psi_{p,\pm}, z_{p\,pm}) 
\end{equation} 
satisfies system \eqref{system:limit}, below. 
\end{thm}
\begin{proof}
This is done in section \ref{sec:convergence_total}. In particular, the strong convergence can be found in
\eqref{cnvg:strong}, \eqref{cnvg:strong-1}, \eqref{cnvg:strong-2}, and \eqref{cnvg:strong-3}, respectively. 
\end{proof}
\begin{remark}
    In this paper, we have not explored the well-posedness, in particular, the uniqueness, of solutions to the limit system \eqref{system:limit}. For this reason, we only have the subsequence convergence in Theorem \ref{thm:convergence}. However, if one manages to show the well-posedness of solutions to system \eqref{system:limit}, the convergence should be of the whole sequence of $ \varepsilon \rightarrow 0^+ $.
\end{remark}

\bigskip

The rest of this paper is organized as follows. In section \ref{sec:preliminaries}, some preliminaries will be provided, including the notations and a boundary-to-domain extension (lifting) Lemma. The classical quasi-geostrophic approximation with only slow waves, i.e., well-prepared initial data, will be reviewed in section \ref{sec:geostrophic}. The key linear slow-fast waves structure will be discussed in section \ref{sec:waves_linear}. Section \ref{sec:uniform_est} is dedicated to the proof of Theorem \ref{thm:uniform}. This paper will finish with the proof of Theorem \ref{thm:convergence} in section \ref{sec:convergence_total}.

\section{Preliminaries}\label{sec:preliminaries}

\subsection{Notations and an extension Lemma}\label{sec:notatons}

In this paper, we have been and will be using 
\begin{equation}
    \biggl( \begin{array}{c} X_1 \\ X_2 \end{array} \biggr)^\perp = \biggl( \begin{array}{c} - X_2 \\ X_1 \end{array} \biggr)
\end{equation}
to denote the rotation of a two-dimensional vector. $ \dvh $ and $ \crlh $ represent the horizontal divergence and curl operators, respectively. Then for any two-dimensional vector field $ X = (X_1, X_2)^\top $, one has
\begin{equation}
    \dvh X^\perp = - \crlh X \qquad \text{and} \qquad \crlh X^\perp = \dvh X.
\end{equation}

For any functions $ A $ and $ B $, the $\mathcal X$ norms are written as
\begin{equation}
    \norm{A, B}{\mathcal X} = \norm{A}{\mathcal X} + \norm{B}{\mathcal X}. 
\end{equation}

We will use $ \Delta_D^{-1} $ to represent the inverse Laplacian subject to the Dirichlet boundary condition at $ z = 0,h $ and the periodic boundary condition horizontally, i.e., 
\begin{equation}
    \Delta \Delta_D^{-1} A = A \qquad \text{with} \quad (\Delta_D^{-1}A)\vert_{z=0,h} = 0.
\end{equation}
Therefore, the definition implies
\begin{equation}\label{id:inverse-Delta}
    \Delta \Delta_D^{-1} = \mathrm{Id}.
\end{equation}
However, observe that \begin{equation}\label{operator-delta-d}
    \Delta_D^{-1} \Delta \neq \mathrm{Id},
\end{equation}
which plays an important role in the proof of short time stability of analytic Prandtl boundary layer \cite{Maekawa2014,nguyenInviscidLimitNavier2018}.

Moreover, $ \Deltah^{-1} $ is the inverse Laplacian in the horizontal variable with zero mean value. Therefore, one has that
\begin{equation}\label{id:inverse-Deltah}
    \Deltah^{-1} \Deltah A = 
    A - \int_{\mathbb T^2} A \,dxdy.
\end{equation}

We will need the following extension (lifting) Lemma:
\begin{lm}\label{lm:extension}
There exists a bi-linear extension operator 
\begin{equation}\label{def:extension}
    \mathrm E_b: \mathcal D'(\mathbb T^2) \times \mathcal D'(\mathbb T^2) \mapsto \mathcal D'(\Omega), 
\end{equation}
such that
for any $ A, B \in H^{s-\frac{1}{2}}(\mathbb T^2) $, $ \mathrm E_b(A,B) \in H^s(\Omega) $ satisfying
\begin{equation}\label{def:extension-to-interial}
    \norm{\mathrm E_b(A,B)}{H^s(\Omega)} \leq C_s \norm{A,B}{H^{s-1/2}(\mathbb T^2)}, 
\end{equation}
and
\begin{equation}\label{def:extension-1}
    \mathrm E_b(A,B)\vert_{z=0} = A \qquad \text{and} \qquad \mathrm E_b(A,B)\vert_{z=h} = B. 
\end{equation}
Moreover, the following property holds:
\begin{equation}
    \dt \mathrm E_b(A,B) = \mathrm E_b(\dt A,\dt B).
\end{equation}
\end{lm}
\begin{proof}
Let $ \chi_0:[0,h] \rightarrow [0,1] $ be a $ C^\infty([0,h]) $ monotonic function
such that
\begin{equation}
    \chi_0(z) = \begin{cases}
    1 & \text{in} \quad z \in [0,h/4), \\
    0 & \text{in} \quad z \in (3h/4,h].
    \end{cases}
\end{equation}
Denote by, $ \vec x_h = (x,y)^\top \in \mathbb T^2 $, for $ A, B \in \mathcal D'(\mathbb T^2) $, 
\begin{equation}
    A(x,y) = \sum_{\vec k \in \mathbb Z^2}A_k e^{i 2\pi \vec k \cdot \vec x_h}, \qquad \text{and} \qquad B(x,y) = \sum_{\vec k \in \mathbb Z^2}B_k e^{i 2\pi \vec k \cdot \vec x_h}.
\end{equation}
For $ z \in \lbrack 0, h \rbrack $, we define
\begin{equation}
    \begin{aligned}
        \mathrm E_b(A,B) = & \sum_{\vec k \in \mathbb Z^2} A_k e^{i 2\pi \vec k \cdot \vec x_h} e^{-|\vec k|z}  \chi_0(z) \\
        & \quad + \sum_{\vec k \in \mathbb Z^2} B_k e^{i 2\pi \vec k \cdot \vec x_h} e^{-|\vec k|(h-z)}( 1 - \chi_0(z)).
    \end{aligned}
\end{equation}
Then it is easy to verify that $ \mathrm E_b(A,B) $ satisfies the properties in the Lemma. This finishes the proof. 
\end{proof}

\subsection{Classical quasi-geostrophic approximation and the potential vorticity formulation for inviscid flows}\label{sec:geostrophic}

In this section, we review the formal quasi-geostrophic approximation with only slow waves of system \eqref{sys:ns_rt_by}, i.e., with well-prepared initial data. This is done by first introducing the formal asymptotic expansion ansatz
\begin{equation}\label{def:asym-epsn}
    \psi(x,y,z,t) : = \psi^0(x,y,z,t) + \varepsilon \psi^1(x,y,z,t)
\end{equation}
for $ \psi \in \lbrace v, w, p, \theta \rbrace $. Then, after substituting \eqref{def:asym-epsn} in system \eqref{sys:ns_rt_by} and matching the $ \mathcal{O}(\varepsilon^{-1}) $ and $ \mathcal O(1) $ terms, one has
\begin{gather}
    (v^0)^\perp + \nablah p^0 = 0, \label{eq:gtph-v}\\
    \dz p^0 - \theta^0 = 0, \label{eq:gtph-w}\\
    w^0 = 0, \label{eq:gtph-tm} \\
    \dt v^0 + v^0 \cdot \nablah v^0 + w^0 \dz v^0 + (v^1)^\perp + \nablah p^1 = 0, \label{eq:sl-v} \\
    \dt w^0 + v^0 \cdot \nablah w^0 + w^0 \dz w^0 + \dz p^1 - \theta^1 = 0, \label{eq:sl-w} \\
    \dt \theta^0 + v^0 \cdot \nablah \theta^0 + w^0 \dz \theta^0 + w^1 = 0, \label{eq:sl-tm} \\
    \dvh v^0 + \dz w^0 = 0, \label{eq:sl-dv} \\
    \shortintertext{and}
    w^0\vert_{z=0,h} = 0. \label{bc:sl}
\end{gather}
In addition, the $ \mathcal O(\varepsilon) $ terms of \eqref{eq:ns-dv} and \eqref{bc:ns} yield
\begin{gather}
    \label{eq:sl-dv-1}
    \dvh v^1 + \dz w^1 = 0,\\
    \shortintertext{and}
    \label{bc:sl-1}
    w^1\vert_{z=0,h} = 0.
\end{gather}

\bigskip

Following \cite{bourgeoisValidityQuasigeostrophicModel1994, embidAveragingFastGravity1996}, we introduce the potential vorticity formulation. Indeed, from \eqref{eq:gtph-v} and \eqref{eq:gtph-w}, it follows that
\begin{equation}\label{def:potential-vorticity}
    \Delta p^0 = (\Deltah + \partial_{zz}) p^0 = \crlh v^0 + \partial_z \theta^0.
\end{equation}
In particular, the quantity on the right hand side of \eqref{def:potential-vorticity} is referred to as the potential vorticity in the literature, and $ p^0 $ is the corresponding steam function. In fact, this terminology is justified by observing that the potential vorticity is transported (see \eqref{eq:ptl-vt}, below). 
After applying $ \crlh $ to \eqref{eq:sl-v}, $ \dz $ to \eqref{eq:sl-tm}, and summing up the resulting equations, one arrives at Ertel's conservation (transport) of the potential vorticity, i.e., 
\begin{equation}\label{eq:ptl-vt}
    \dt \Delta p^0 + v^0 \cdot \nablah \Delta p^0  = 0,
\end{equation}
where we have applied the fact, thanks to \eqref{eq:gtph-v}, \eqref{eq:gtph-w}, \eqref{eq:gtph-tm}, and \eqref{eq:sl-dv}, that 
\begin{equation}\label{eq:sl-dv-2}
    \dz v^0 \cdot \nablah \theta^0 = 0, \qquad
    w^0 = 0, \qquad \text{and} \qquad \dvh v^0 = 0.
\end{equation}
In addition, thanks to \eqref{eq:gtph-w}, \eqref{eq:sl-tm}, and \eqref{bc:sl-1}, one can show that
\begin{equation}\label{eq:sl-bc}
    \dt (\dz p^0\vert_{z=0,h}) + v^0 \vert_{z=0,h} \cdot \nablah (\dz p^0\vert_{z=0,h}) = 0.
\end{equation}

\bigskip

The system formed by \eqref{eq:gtph-v}, \eqref{eq:gtph-w}, \eqref{eq:ptl-vt}, and \eqref{eq:sl-bc} is the well-known potential vorticity formulation of the classical quasi-geostrophic approximation. In particular, \eqref{eq:sl-bc} describes the evolution of `boundary conditions' for the stream function $ p^0 $, i.e., $ \dz p^0\vert_{z=0,h} $, which is used to invert the Laplacian in $ v^0 = \nablah^\perp p^0 = \nablah^\perp \Delta_N^{-1} (\Delta p^0) $, where $ \Delta_N^{-1} $ here is the inverse Laplacian with Neumann type boundary condition at $ z = 0,h $ and periodic boundary condition horizontally. Observe from \eqref{eq:sl-bc} that if $ \dz p^0\vert_{z=0,h} = 0 $ initially, it remains zero. This is one of the underlying observation behind the well-prepared initial data in \cite{bourgeoisValidityQuasigeostrophicModel1994}. In addition, observe that $ \Delta_N^{-1} $ is unique up to a constant, which, without loss of generality, can be taken to be zero, justifying the notation of inverse.

\subsection{The slow--fast waves structure: Linear analysis}\label{sec:waves_linear}

Our goal in this section is to investigate the linear slow-fast waves structure of system \eqref{sys:ns_rt_by}. This will guide us to obtain uniform-in-$\varepsilon$ estimates as well as nonlinear waves interaction analysis in the next sections. Without loss of generality, we write $ (v_l, w_l, \theta_l) $ and $ p_l $, i.e., the linear variables, and the linear system associated with system \eqref{sys:ns_rt_by} as follows:
\begin{subequations}\label{sys:ns_rt_by_l}
\begin{align}
    \dt v_l && + \dfrac{1}{\varepsilon} v_l^\perp && + \dfrac{\nablah p_l}{\varepsilon} && && = 0, \label{eq:ns_v_l}\\
    \dt w_l && && + \dfrac{\dz p_l}{\varepsilon} && - \dfrac{\theta_l}{\varepsilon} && = 0, \label{eq:ns_w_l} \\
    \dt \theta_l && && && + \dfrac{w_l}{\varepsilon} && = 0, \label{eq:ns_tm_l} \\
    \dvh v_l + \dz w_l && && && && = 0, \label{eq:ns_dv_l}
\end{align}
with
\begin{equation}
    w_l\vert_{z=0,h} = 0 \qquad \text{i.e., ~~ impermeable boundary condition}, \label{bc:ns_l}
\end{equation}
and periodic boundary condition horizontally. 
\end{subequations}

\bigskip

The linear version of Ertel's conservation (transport) of the potential vorticity $(\dz \theta_l + \crlh v_l)$
and the corresponding stream function $ p_l $
read, thanks to \eqref{eq:ns_v_l}, \eqref{eq:ns_dv_l}, and \eqref{bc:ns_l},
\begin{subequations}\label{sys:linear-fast-slow-dcmp}
\begin{equation}\label{eq:potential_vorticity_l}
   \Delta_h p_l + \partial_{zz} p_l = \dz \theta_l + \crlh v_l, \quad \dt(\Delta_h p_l + \partial_{zz} p_l) = \dt (\dz \theta_l + \crlh v_l) = 0.
\end{equation}
Meanwhile, taking the trace of \eqref{eq:ns_tm_l} to the channel boundary yields
\begin{equation}\label{eq:theta-bc-l}
    \dt \theta_l\vert_{z=0,h} = 0.
\end{equation}
On the other hand, one can verify that
\begin{equation}\label{eq:generalized_pv_l}
    \dt (\nablah^\perp \theta_l + \nablah w_l - \dz v_l) + \dfrac{1}{\varepsilon}(\nablah^\perp \theta_l + \nablah w_l - \dz v_l)^\perp = 0.
\end{equation}
Last but not least, integrating \eqref{eq:ns_v_l} in the horizontal variables yields
\begin{equation}\label{eq:mean_v_l}
    \dt \int_{\mathbb T^2} v_l(x,y,z)\,dxdy + \dfrac{1}{\varepsilon}\biggl(\int_{\mathbb T^2} v_l(x,y,z)\,dxdy\biggr)^\perp = 0.
\end{equation}

\end{subequations}

Moreover, observe that \eqref{eq:ns_w_l} and \eqref{eq:ns_tm_l} imply 
\begin{equation}\label{eq:ns_bc_p_l}
    \dt (\dz p_l\vert_{z=0,h}) = 0.
\end{equation}

\bigskip

Equations \eqref{eq:potential_vorticity_l} and \eqref{eq:generalized_pv_l} form the linear full generalized potential vorticity equations. 
A few remarks about this linear structure are in order:
\begin{itemize}
    \item While system \eqref{sys:ns_rt_by_l} is stable with respect to the $ L^2 $ norm, i.e., one can get uniform-in-$\varepsilon $ $ L^2 $ estimate by taking the $ L^2 $-inner product of \eqref{eq:ns_v_l}, \eqref{eq:ns_w_l}, and \eqref{eq:ns_tm_l} with respect to $ v_l $, $ w_l $, and $ \theta_l $, the same can not be said about the $ H^s $ estimate for $ s \geq 1 $. This is due to the absence of boundary condition for the higher order derivatives of $ p_l $ and $ w_l $. For this reason, only in the case of periodic spatial domains (e.g., \cite{embidAveragingFastGravity1996}), or in the case with well-prepared initial data and $ \partial_z p_l\vert_{z=0,h} = 0 $ (e.g., \cite{bourgeoisValidityQuasigeostrophicModel1994}; see \eqref{eq:ns_bc_p_l}), one can verify the uniform $H^s$ estimates and the asymptotic limit  as $ \varepsilon \rightarrow 0^+ $; 
    
    \item On the other hand, \eqref{eq:potential_vorticity_l}, \eqref{eq:generalized_pv_l}, and \eqref{eq:mean_v_l} completely eliminate $ p_l $, and in particular, the underlying quantities in this system are stable with respect to any spatial derivatives. Therefore, one can get uniform-in-$\varepsilon$ $ H^s $ estimates without any restriction for these quantities;
     
    \item To be more precise, the estimates of the horizontal derivatives can be achived from \eqref{sys:ns_rt_by_l}. Then from \eqref{eq:potential_vorticity_l}, \eqref{eq:generalized_pv_l}, and \eqref{eq:ns_dv_l}, one can derive the estimates of $ \dz \theta_l $, $ \dz v_l $, and $ \dz w_l $, respectively, in terms of the horizontal derivatives. Bootstrap arguments will lead to $ H^s $ estimates;
    
    \item One can regard \eqref{eq:potential_vorticity_l} and \eqref{eq:theta-bc-l} as the equations of the slow waves (dynamics), and \eqref{eq:generalized_pv_l} and \eqref{eq:mean_v_l} as the equations of the fast waves (dynamics). That is, one is able to separate the slow and fast state variables;
    
    \item From \eqref{eq:ns_tm_l} and \eqref{eq:generalized_pv_l}, one can conclude that as $ \varepsilon \rightarrow 0 $, $ w_l, \nablah^\perp \theta_l - \dz v_l \rightharpoonup  0 $, weakly in the sense of distribution. This is consistent with 
    \eqref{eq:gtph-v}, \eqref{eq:gtph-w}, and
    \eqref{eq:gtph-tm}. 
\end{itemize}

\bigskip

Now we shall write down the slow-fast waves of linear system \eqref{sys:ns_rt_by_l}. Denote by
\begin{align}
    \label{def:pv_l} \Phi_l(x,y,z,t) := & \dz \theta_l + \crlh v_l  \qquad (\text{the potential vorticity}), \\
    \label{def:g_pv_l} \Psi_l(x,y,z,t) := & \nablah^\perp \theta_l + \nablah w_l - \dz v_l, \\
    \label{def:bc_th_l-0} H_{l,0}(x,y,t) := & \theta_l\vert_{z=0}, \\
    \label{def:bc_th_l-h} H_{l,h}(x,y,t) := & \theta_l\vert_{z=h}, \\
    \intertext{and}
    \label{def:bc_v_l} Z_l(z,t) := & \int_{\mathbb T^2} v_l(x,y,z) \,dxdy.
\end{align}
Correspondingly, let $ \Phi_\mrm{in} $, $\Psi_\mrm{in}$, $H_{0,\mrm{in}}$, $H_{h,\mrm{in}}$, and $Z_\mrm{in}$ be the initial data at $ t = 0 $ for $ \Phi_l $, $ \Psi_l $, $ H_{l,0} $, $ H_{l,h} $, and $ Z_l $, respectively. 
In particular, $ \Phi_l $ and $ \Psi_l $ form the generalized potential vorticity, and are the main ingredient of, and to be explored later in, this work. 
Then it follows from system \eqref{sys:linear-fast-slow-dcmp},
that
\begin{equation}\label{sol:linear-sys}
\begin{aligned}
    \text{linear slow variables:} & \quad \Phi_l(t) \equiv \Phi_\mrm{in}, \quad H_{l,0}(t) \equiv H_{0,\mrm{in}}, \quad H_{l,h}(t) \equiv H_{h,\mrm{in}}, \\
    \text{linear fast variables:} & \quad
    \Psi_l(t) = e^{it/\varepsilon} \dfrac{\Psi_\mrm{in} + i \Psi_\mrm{in}^\perp}{2} + e^{-it/\varepsilon} \dfrac{\Psi_\mrm{in} - i \Psi_\mrm{in}^\perp}{2},\\ 
    \text{and} & \quad Z_l(t) = e^{it/\varepsilon} \dfrac{Z_\mrm{in} + i Z_\mrm{in}^\perp}{2} + e^{-it/\varepsilon} \dfrac{Z_\mrm{in} - i Z_\mrm{in}^\perp}{2}.
\end{aligned}
\end{equation}

We claim that $ (\Phi_l, \Psi_l, H_{l,0}, H_{l,h}, Z_l) $ as in \eqref{sol:linear-sys} provide complete information on the solutions of system \eqref{sys:ns_rt_by_l}. This can be seen by writing $ (v_l, w_l, \theta_l) $ in terms of $ (\Phi_l, \Psi_l, H_{l,0}, H_{l,h}, Z_l) $. First, taking $ \dvh $ and $ \crlh $ to \eqref{def:g_pv_l} yields that, respectively, thanks to \eqref{eq:ns_dv_l} and \eqref{def:pv_l},
\begin{gather}
    \Deltah w_l + \partial_{zz} w_l = \dvh \Psi_l \label{eq:laplace_w_l} \\
    \intertext{and}
    \Deltah \theta_l + \partial_{zz} \theta_l = \partial_{z} \Phi_l + \crlh \Psi_l \qquad \text{or, equivalently} \nonumber\\
    \Delta ( \theta_l - \mathrm E_{b}(H_{l,0},H_{l,h})) = \crlh \Psi_l + \dz \Phi_l - \Delta \mathrm E_{b}(H_{l,0},H_{l,h}). \label{eq:laplace_tm_l}
\end{gather}
Note that, thanks to  \eqref{def:extension-1}, \eqref{bc:ns_l}, \eqref{def:bc_th_l-0}, and \eqref{def:bc_th_l-h},
\begin{equation*}
    w_l\vert_{z=0,h} = 0 \qquad\text{and} \qquad (\theta_l - \mathrm E_{b}(H_{l,0},H_{l,h}))\vert_{z=0,h} = 0.
\end{equation*}
Therefore, let $ \Delta^{-1}_D $ be the three-dimensional inverse Laplacian with Dirichlet boundary condition on $ \lbrace z = 0,h \rbrace $ and periodic boundary condition in the horizontal directions. From \eqref{eq:laplace_w_l} and \eqref{eq:laplace_tm_l}, one has
\begin{align}
    w_l = & \Delta_D^{-1} \dvh \Psi_l \label{sol:w_l} \\
    \intertext{and}
    \theta_l = & \mathrm E_{b}(H_{l,0},H_{l,h}) + \Delta_D^{-1}(\crlh \Psi_l + \dz  \Phi_l - \Delta \mathrm E_{b}(H_{l,0},H_{l,h})). \label{sol:tm_l}
\end{align}
To calculate $ v_l $, let $ \Deltah^{-1} $ be the two-dimensional inverse Laplace with zero horizontal mean value. Then, thanks to \eqref{eq:ns_dv_l} and \eqref{def:pv_l}, one has 
\begin{equation}\label{eq:dv-curl-v-l}
    \dvh v_l = - \dz w_l \qquad \text{and} \qquad \crlh v_l =  \Phi_l - \dz \theta_l,
\end{equation}
and, therefore, it follows that
\begin{equation*}
    v_l = Z_l + \nablah \Deltah^{-1} \dvh v_l + \nablah^\perp \Deltah^{-1}\crlh v_l,
\end{equation*}
or, after substituting \eqref{eq:dv-curl-v-l}, \eqref{sol:w_l}, and \eqref{sol:tm_l} in the above expression, one has
\begin{equation}\label{sol:v_l}
\begin{aligned}
    v_l = & Z_l - \nablah \Deltah^{-1} \dz(\Delta_D^{-1}\dvh \Psi_l)\\
     & + \nablah^\perp \Deltah^{-1} [\Phi_l - \dz \mathrm E_{b}(H_{l,0},H_{l,h}) \\
     & \qquad \qquad - \dz \Delta_D^{-1}(\crlh \Psi_l + \dz  \Phi_l - \Delta \mathrm E_{b}(H_{l,0},H_{l,h}))].
\end{aligned}
\end{equation}
We remind the reader that $ (\Phi_l, \Psi_l, H_{l,0},H_{l,h}, Z_l) $ are as in \eqref{sol:linear-sys}, with $ (\Psi_l, Z_l) $ being fast state variables and $ (\Phi_l, H_{l,0},H_{l,h}) $ slow state variables. Therefore, one can decompose $ v_l,w_l, \theta_l $ in terms of slow and fast waves in an unambiguous fashion.

\section{Uniform-in-$ \varepsilon $ estimates of the Euler equations with fast Rossby and gravity waves}\label{sec:uniform_est}

In this and the following sections, we will proceed to the nonlinear analysis. In particular, we focus in this section on the uniform-in-$\varepsilon$ estimates for system \eqref{sys:ns_rt_by} in this section. Inspired by the discussion in section \ref{sec:waves_linear}, we define
\begin{align}
    \label{def:pv} \Phi(x,y,z,t): = & \dz \theta + \crlh v, \\
    \label{def:g_pv}  \Psi(x,y,z,t): = & \nablah^\perp \theta + \nablah w - \dz v, \\
    \label{def:bc_th_0} H_0(x,y,t) := & \theta\vert_{z=0}, \\
    \label{def:bc_th_h}
    H_h(x,y,t) := & \theta\vert_{z=h}, \\
    \shortintertext{and}
    \label{def:bc_v} Z(z,t) := & \int_{\mathbb{T}^2} v(x,y,z,t) \,dxdy.
\end{align}
Recall that $ \Phi $ and $ \Psi $ form the {\bf generalized potential vorticity}. 
From \eqref{eq:ns_v}, \eqref{eq:ns_w}, \eqref{eq:ns_tm}, and \eqref{eq:ns-dv}, one can write down the following equations
\begin{gather}
    \label{eq:curlh_v} \begin{gathered} \dt \crlh v + v \cdot \nablah \crlh v + w \dz \crlh v \\ + \crlh v \cdot \dvh v + \dz v \cdot \nablah^\perp w - \dfrac{\dz w}{\varepsilon} = 0,\end{gathered} \\
    \label{eq:dz_v} \begin{gathered} \dt \dz v + v \cdot \dz v + w \dz \dz v + \dfrac{\dz v^\perp}{\varepsilon} + \dfrac{\nablah \dz p}{\varepsilon} \\ + \dz v \cdot \nablah v + \dz w \dz v  = 0, \end{gathered} \\
    \label{eq:nablah_w} \begin{gathered} \dt \nablah w + v \cdot \nablah \nablah w + w \dz \nablah w + \dfrac{\nablah \dz p}{\varepsilon} - \dfrac{\nablah \theta}{\varepsilon} \\ + (\nablah v)^\top \nablah w + \dz w \nablah w = 0,\end{gathered} \\
    \label{eq:nablah_theta} \begin{gathered} \dt \nablah \theta + v \cdot \nablah \nablah \theta + w \dz \nablah \theta + \dfrac{\nablah w}{\varepsilon} \\ + (\nablah v)^\top \nablah \theta + \dz \theta \nablah w = 0,
    \end{gathered}\\
    \label{eq:dz_theta} \begin{gathered} \dt \dz \theta + v \cdot \nablah \dz \theta + w \dz \dz \theta + \dfrac{\dz w}{\varepsilon} \\ + \dz v \cdot \nablah \theta + \dz w \dz \theta  = 0.\end{gathered}
\end{gather}
Consequently, one has, from system \eqref{sys:ns_rt_by}, that 
\begin{subequations}\label{sys:ns_rt_by_s_l}
    \begin{gather}
        \label{eq:pv} \begin{gathered} 
            \dt \Phi + v \cdot \nablah  \Phi + w \dz  \Phi 
            + N_1 = 0,
        \end{gathered} \\
        \label{eq:g_pv} \begin{gathered}
            \dt \Psi + v \cdot \nablah \Psi + w \dz \Psi + \dfrac{1}{\varepsilon} \Psi^\perp + N_2 = 0,
        \end{gathered}\\
        \label{eq:bc_th_0} \begin{gathered} \dt H_0 + v\vert_{z=0} \cdot \nablah H_0 = 0,\end{gathered} \\
        \label{eq:bc_th_h} \begin{gathered} \dt H_h + v\vert_{z=h} \cdot \nablah H_h = 0,\end{gathered} \\
        \label{eq:bc_v} \dt Z + \dfrac{1}{\varepsilon} Z^\perp + N_3 = 0,
    \end{gather}
    where
    \begin{align}
        \label{def:nonln_pv} N_1:= & \crlh v \cdot \dvh v + \dz v \cdot \nablah^\perp w 
         + \dz v \cdot \nablah \theta + \dz w \dz \theta, \\
        \label{def:nonln_g_pv} N_2 := & ((\nablah v)^\top \nablah \theta)^\perp + \dz \theta \cdot \nablah^\perp w + (\nablah v)^\top \nablah w + \dz w \nablah w \nonumber \\
        & - \dz v \cdot \nablah v - \dz w \dz v, \\
        \label{def:nonln_bc_v} N_3:= & \int_{\mathbb T^2} \dz (w v) \,dxdy.
    \end{align}
\end{subequations}

\bigskip

We continue with the uniform-in-$\varepsilon$ estimates in the following steps: 1. establish estimates for the horizontal derivatives; then 2. establish estimates for the vertical derivatives; finally, 3. close the estimates.

\subsection*{Estimates for the horizontal derivatives}

Let $ \ddh \in \lbrace \partial_x, \partial_y \rbrace $ and $ \alpha \in \lbrace 0,1,2,3 \rbrace $. Applying $ \ddh^\alpha $ to system \eqref{sys:ns_rt_by} leads to
    \begin{gather}
        \label{eq:ddh_v} \begin{gathered}
            \dt \ddh^\alpha v + (v \cdot \nablah + w \dz) \ddh^\alpha v + \dfrac{1}{\varepsilon} \ddh^\alpha v^\perp + \dfrac{\nablah\ddh^\alpha p}{\varepsilon} \\
            + \ddh^\alpha (v\cdot \nablah v + w \dz v) - (v\cdot \nablah + w \dz) \ddh^\alpha v = 0,
        \end{gathered}\\
        \label{eq:ddh_w} \begin{gathered}
            \dt \ddh^\alpha w + (v \cdot \nablah + w \dz) \ddh^\alpha w + \dfrac{\dz \ddh^\alpha p}{\varepsilon} - \dfrac{\ddh^\alpha \theta}{\varepsilon} \\
            + \ddh^\alpha (v\cdot \nablah w + w \dz w) - (v \cdot \nablah + w \dz) \ddh^\alpha w = 0,
        \end{gathered} \\
        \label{eq:ddh_tm} \begin{gathered}
            \dt \ddh^\alpha \theta + (v \cdot \nablah + w \dz ) \ddh^\alpha \theta + \dfrac{\ddh^\alpha w}{\varepsilon} \\ + \ddh^\alpha (v\cdot\nablah \theta + w \dz \theta) - (v \cdot \nablah + w \dz ) \ddh^\alpha \theta = 0,
        \end{gathered}\\
        \dvh \ddh^\alpha v + \dz \ddh^\alpha w = 0, \qquad \ddh w\vert_{z=0,h} = 0.
    \end{gather}
Taking the $ L^2 $-inner product of \eqref{eq:ddh_v}--\eqref{eq:ddh_tm}  with $ 2\ddh^\alpha v, 2\ddh^\alpha w, 2\ddh^\alpha \theta $, respectively, applying integration by parts, and summing up the resultants lead to
\begin{equation}\label{eq:energy_hh}
    \begin{gathered}
        \dfrac{d}{dt} \norm{\ddh^\alpha v, \ddh^\alpha w, \ddh^\alpha \theta}{\Lnorm{2}}^2 \\
        = 
        - 2\int \lbrack \ddh^\alpha (v\cdot \nablah v + w \dz v) - (v\cdot \nablah + w \dz) \ddh^\alpha v \rbrack \cdot \ddh^\alpha v \idx \\
        - 2\int \lbrack \ddh^\alpha (v\cdot \nablah w + w \dz w) - (v \cdot \nablah + w \dz) \ddh^\alpha w \rbrack \times \ddh^\alpha w \idx \\
        - 2\int \lbrack \ddh^\alpha (v\cdot\nablah \theta + w \dz \theta) - (v \cdot \nablah + w \dz ) \ddh^\alpha \theta \rbrack \times \ddh^\alpha \theta\\
        \leq C \norm{v, w, \theta}{\Hnorm{2}}^{1/2}\times\norm{v, w, \theta}{\Hnorm{3}}^{5/2},
    \end{gathered}
\end{equation}
for some generic constant $ C \in (0,\infty) $, 
where in the last inequality we have applied the H\"older inequality, the Gagliardo-Nirenberg inequality, and the Sobolev embedding inequality.

\subsection*{Estimates for the vertical derivatives}
As before, let $ \partial \in \lbrace \partial_x, \partial_y, \partial_z \rbrace $ and $ \beta \in \lbrace 0,1,2 \rbrace $. Applying $ \partial^\beta $ to equations \eqref{eq:pv} and \eqref{eq:g_pv} leads to
\begin{gather}
    \label{eq:dd_pv} \begin{gathered}
        \dt \partial^\beta \Phi + (v \cdot \nablah + w \dz ) \partial^\beta  \Phi  + \partial^\beta N_1\\
        + \partial^\beta (v \cdot \nablah \Phi + w \dz \Phi ) - (v \cdot \nablah + w \dz ) \partial^\beta \Phi = 0,
    \end{gathered} \\
    \label{eq:dd_g_pv} \begin{gathered}
        \dt \partial^\beta \Psi + (v\cdot \nablah + w \dz ) \partial^\beta \Psi + \dfrac{1}{\varepsilon} \partial^\beta \Psi^\perp + \partial^\beta N_2 \\
        + \partial^\beta (v \cdot \nablah \Psi + w \dz \Psi) - (v \cdot \nablah + w \dz )\partial^\beta \Psi = 0.
    \end{gathered}
\end{gather}
Taking the $ L^2 $-inner product of \eqref{eq:dd_pv} and \eqref{eq:dd_g_pv} with $ 2\partial^\beta \Phi $ and $ 2\partial^\beta \Psi $, respectively, applying integration by parts, and summing up the resultants lead to
\begin{equation}\label{eq:energy_dd}
    \begin{gathered}
        \dfrac{d}{dt} \norm{\partial^\beta  \Phi, \partial^\beta \Psi}{\Lnorm{2}}^2  = 
        - 2 \int \bigl( \partial^\beta N_1 \cdot \partial^\beta  \Phi + \partial^\beta N_2 \cdot \partial^\beta \Psi \bigr) \idx \\
        - 2 \int \lbrack \partial^\beta (v \cdot \nablah  \Phi + w \dz  \Phi ) - (v \cdot \nablah + w \dz ) \partial^\beta  \Phi \rbrack \cdot \partial^\beta  \Phi \idx \\
        - 2 \int \lbrack \partial^\beta (v \cdot \nablah \Psi + w \dz \Psi) - (v \cdot \nablah + w \dz )\partial^\beta \Psi \rbrack \cdot \partial^\beta \Psi \idx \\
        \leq  C \norm{v,w,\theta}{\Hnorm{3}}^2 \norm{\Phi, \Psi}{\Hnorm{2}} + C \norm{v,w,\theta}{\Hnorm{3}} \norm{\Phi, \Psi}{\Hnorm{2}}^2,
    \end{gathered}
\end{equation}
for some absolute constant $ C \in (0,\infty) $, where in the last inequality we have applied the H\"older inequality, the Gagliardo-Nirenberg inequality, and the Sobolev embedding inequality.

\subsection*{Closing the estimates} 
Define the total ``energy'' functional by 
\begin{equation}\label{def:energy_funtl}
    \mathfrak E := \norm{\Phi, \Psi}{\Hnorm{2}}^2 + \sum_{\mathclap{\substack{\ddh \in \lbrace \partial_x, \partial_y \rbrace, \\ \alpha\in \lbrace 0,1,2,3 \rbrace }}} \norm{\ddh^\alpha v, \ddh^\alpha w, \ddh^\alpha \theta}{\Lnorm{2}}^2.
\end{equation}
We observe that
\begin{equation}\label{eqvl:energy_funtl}
    \dfrac{1}{C} \norm{v,w,\theta}{\Hnorm{3}}^2 \leq \mathfrak E \leq C \norm{v,w,\theta}{\Hnorm{3}}^2,
\end{equation}
for some generic constant $ C \in (0,\infty) $. Indeed, the right-hand side inequality in \eqref{eqvl:energy_funtl} follows directly from the definition of $ \Phi $ and $ \Psi $ in \eqref{def:pv} and \eqref{def:g_pv}. To show the left-hand side inequality, notice that
\begin{gather*}
    \partial_z v = - \Psi + \nablah^\perp \theta + \nablah w, \qquad 
    \dz \theta = \Phi - \crlh v, \\
    \shortintertext{and}
    \dz w = - \dvh v. 
\end{gather*}
Thus,
\begin{align*}
    \sum_{\alpha\in \lbrace 0,1,2 \rbrace}\norm{\ddh^\alpha \dz v, \ddh^\alpha \dz w, \ddh^\alpha \dz \theta}{\Lnorm{2}} \leq C \mathfrak E.
\end{align*}
Similarly, following a bootstrap argument on the derivatives implies the left-hand side part of \eqref{eqvl:energy_funtl}.

Consequently, \eqref{eq:energy_hh} and \eqref{eq:energy_dd} yield
\begin{equation}\label{ineq:uniform_est}
    \dfrac{d}{dt} \mathfrak E \leq C \mathfrak E^{3/2},
\end{equation}
for some generic constant $ C \in (0,\infty) $. In particular, from \eqref{ineq:uniform_est} and \eqref{eqvl:energy_funtl}, one concludes that there exists $ T \in (0,\infty) $, depending only on the initial data and independent of $ \varepsilon $, such that
\begin{equation}\label{est:uniform}
    \sup_{0\leq t\leq T} \norm{v(t),w(t),\theta(t)}{\Hnorm{3}}^2 \leq C \sup_{0\leq t\leq T} \mathfrak E(t) \leq 2 C^2 \norm{v_\mrm{in},w_\mrm{in},\theta_\mrm{in}}{\Hnorm{3}}^2,
\end{equation}
for the same constant $ C $ as in \eqref{eqvl:energy_funtl}. This finishes the proof of Theorem \ref{thm:uniform}.

\section{Convergence theory}\label{sec:convergence_total}

\subsection{Convergence theory: Part 1,  compactness}\label{sec:convergence_1}

What is left is to establish the convergence of the solutions to system \eqref{sys:ns_rt_by} as $ \varepsilon \rightarrow 0^+ $, which we will do in two steps. In this subsection, we will conclude the weak and strong compactness, thanks to the uniform estimate \eqref{est:uniform}. In the next subsection, we will deal with the convergence of the nonlinearities.

In the rest of this paper, we denote by $ T \in (0,\infty) $ the uniform-in-$\varepsilon$ existence time established in section \ref{sec:uniform_est} at \eqref{est:uniform}. $ C_\mrm{in}\in (0,\infty) $ will denote a constant that is independent of $ \varepsilon $,
different from line to line, depending only on the initial data. With such notations, thanks to \eqref{est:uniform}, by virtue of the definitions of $ \Phi $, $ \Psi $, $ H_0 $, $ H_h $, and $ Z $ in \eqref{def:pv}--\eqref{def:bc_v}, respectively, we have
\begin{equation}\label{uni-est:001}
    \sup_{0\leq t \leq T} \bigl( \norm{\Phi(t), \Psi(t)}{\Hnorm{2}} + \norm{H_0(t), H_h(t)}{H^{5/2}(\mathbb T^2)} + \norm{Z(t)}{\Hnorm{3}} \bigr) \leq C_\mrm{in}.
\end{equation}
Similarly, from \eqref{def:nonln_pv}--\eqref{def:nonln_bc_v}, it follows that
\begin{equation}\label{uni-est:001-1}
    \sup_{0\leq t \leq T} \bigl( \norm{N_1,N_2,N_3}{\Hnorm{2}} \bigr) \leq C_\mrm{in}.
\end{equation}

From \eqref{eq:pv}--\eqref{eq:bc_v}, one has, thanks to \eqref{est:uniform}, \eqref{uni-est:001}, and \eqref{uni-est:001-1}, that
\begin{equation}\label{uni-est:003}
    \begin{gathered}
    \sup_{0\leq t \leq T}\bigl( \norm{\dt \Phi(t), \varepsilon \dt \Psi(t)}{\Hnorm{1}}+ \norm{\dt H_0(t), \dt H_h(t)}{H^{3/2}(\mathbb T^2)} \\
    + \norm{\varepsilon \dt Z(t)}{\Hnorm{2}}\bigr) \leq C_\mrm{in}.
    \end{gathered}
\end{equation}

Consequently, by virtue of the Aubin compactness theorem \cite[Theorem 2.1]{Temam}, there exist
\begin{equation}\label{def:limit-variables}
\begin{gathered}
    \Phi_p, \Psi_p \in L^\infty(0,T; H^2(\Omega)), \qquad H_{p,0},  H_{p,h} \in L^\infty(0,T; H^{5/2}(\mathbb T^2)), \\ \text{and} \qquad 
    Z_p, v_p, w_p, \theta_p \in L^\infty(0,T;H^3(\Omega)),
\end{gathered}
\end{equation}
with
\begin{equation}\label{def:limit-variables-2}
    \dt \Phi_p \in L^\infty(0,T;H^1(\Omega), \qquad \dt H_{p,0},  \dt H_{p,h} \in L^\infty(0,T;H^{3/2}(\mathbb T^2),
\end{equation}
such that there exists a subsequence of $ \varepsilon $ that as $ \varepsilon \rightarrow 0^+ $,
\begin{align}
    \label{cnvg:weak-1}\Phi, \Psi &\stackrel{\ast}{\rightharpoonup} \Phi_p,  \Psi_p&  \text{weak-$\ast$ in} & \qquad L^\infty(0,T;H^2(\Omega)),\\
    \label{cnvg:weak-1-2} H_0, H_h &\stackrel{\ast}{\rightharpoonup} H_{p,0},  H_{p,h} &  \text{weak-$\ast$ in} & \qquad L^\infty(0,T;H^{5/2}(\mathbb T^2)), \\
    \label{cnvg:weak-2} Z, v, w, \theta &\stackrel{\ast}{\rightharpoonup} Z_p, v_p, w_p, \theta_p &  \text{weak-$\ast$ in} & \qquad L^\infty(0,T;H^3(\Omega)), \\
    \label{cnvg:weak-3} \dt \Phi &\stackrel{\ast}{\rightharpoonup} \dt \Phi_p  &  \text{weak-$\ast$ in} & \qquad L^\infty(0,T;H^1(\Omega))\\
    \label{cnvg:weak-3-2} \dt H_0, \dt H_h &\stackrel{\ast}{\rightharpoonup} \dt H_{p,0}, \dt H_{p,h}  &  \text{weak-$\ast$ in} & \qquad L^\infty(0,T;H^{3/2}(\mathbb T^2))\\
    \shortintertext{and}
    \label{cnvg:strong}\Phi &\rightarrow \Phi_p &  \text{in} & \qquad C([0,T];H^1(\Omega)), \\
    \label{cnvg:strong-1}
    H_0, H_h &\rightarrow H_{p,0}, H_{p,h} & \text{in} & \qquad C([0,T];H^{3/2}(\mathbb T^2))
\end{align}

Furthermore, from \eqref{eq:ns_tm}, \eqref{eq:g_pv} and \eqref{eq:bc_v}, after sending $ \varepsilon \rightarrow 0^+ $, one can verify that $ w_p = \Psi_p = Z_p \equiv 0 $. In fact, after taking the inner product of corresponding equations with $ \varepsilon $ and a test function in $ \mathcal D'((0,T)\times\Omega)$ and passing the limit $ \varepsilon \rightarrow 0^+ $, it is easy to verify that $ w_p = \Psi_p = Z_p \equiv 0 $ in the sense of distribution. Then it follows from the regularity of $ w_p, \Psi_p $, and $ Z_p $ that they are equal to zero. Following similar arguments from the definition, it is easy to show that,
\begin{equation}\label{def:slow-variables}
\begin{gathered}
    w_p = 0, \quad \Phi_p = \dz \theta_p + \crlh v_p, \quad \nablah^\perp \theta_p + \nablah w_p - \dz v_p = 0, \\
    \dvh v_p + \dz w_p = 0, \quad
    H_{p,0} = \theta_p \vert_{z=0}, \quad H_{p,h} = \theta_p\vert_{z=h}, \\ \quad \text{and} \quad \int_{\mathbb{T}^2} v_p(x,y,z) \,dxdy = 0,
\end{gathered}
\end{equation}
or, equivalently, repeating similar calculation as in \eqref{eq:laplace_w_l}--\eqref{sol:v_l}, one has
\begin{equation*}{\tag{\ref{def:slow-variables}'}}\label{def:slow-variables-2}
    \begin{aligned}
        & w_p = 0, && \qquad \theta_p  =  \mathrm E_b(H_{p,0}, H_{p,h}) + \Delta_D^{-1} ( \dz \Phi_p - \Delta \mathrm E_b(H_{p,0}, H_{p,h})),  \\
        & \text{and} && \qquad v_p = \nablah^\perp \Deltah^{-1} [\Phi_p - \dz \mathrm E_b(H_{p,0}, H_{p,h})\\
        & && \qquad\qquad\qquad - \dz \Delta_D^{-1}(\dz \Phi_p - \Delta \mathrm E_b(H_{p,0}, H_{p,h}))].
    \end{aligned}
\end{equation*}

\begin{remark}\label{rk:potential}
    We can perform the following calculation to rewrite $ \theta_p $. 
    Let $ P:= \Deltah^{-1} [\Phi_p - \dz \mathrm E_b(H_{p,0}, H_{p,h}) - \dz \Delta_D^{-1}(\dz \Phi_p - \Delta \mathrm E_b(H_{p,0}, H_{p,h}))] $. Then direct calculation shows that 
    \begin{align*}
        \dz P = &  \Deltah^{-1} [\dz \Phi_p - \partial_{zz} \mathrm E_b(H_{p,0}, H_{p,h})\\
        & \qquad - (\Delta - \Deltah) \Delta_D^{-1}(\dz \Phi_p - \Delta \mathrm E_b(H_{p,0}, H_{p,h}))] \\
        = & \underbrace{\mathrm E_b(H_{p,0}, H_{p,h}) + \Delta_D^{-1} ( \dz \Phi_p - \Delta \mathrm E_b(H_{p,0}, H_{p,h}))}_{=\theta_p}\\
        & \qquad - \underbrace{\int_{\mathbb T^2}\bigl\lbrack \mathrm E_b(H_{p,0}, H_{p,h}) + \Delta_D^{-1} ( \dz \Phi_p - \Delta \mathrm E_b(H_{p,0}, H_{p,h})) \bigr\rbrack \,dxdy}_{=:Q(z)},
    \end{align*}
    where we have applied \eqref{id:inverse-Delta} and \eqref{id:inverse-Deltah}.
    Together with \eqref{def:slow-variables-2}, we have 
    \begin{equation*}
        \theta_p = \dz (P + \int_0^z Q(z')\,dz') \qquad \text{and} \qquad v_p = \nablah^\perp(P + \int_0^z Q(z')\,dz').
    \end{equation*}
    This is consistent with the classical theory of the quasi-geostrophic approximation. See, for instance, \cite{embidAveragingFastGravity1996,bourgeoisValidityQuasigeostrophicModel1994}. 
\end{remark}

\bigskip

Next, to handle the fast waves, i.e., $ \Psi $ and $ Z $, following Schochet's theory \cite{Schochet1994}, from \eqref{eq:g_pv} and \eqref{eq:bc_v}, one has 
\begin{gather}
    \label{eq:fast-wave-1}
    \begin{gathered}
    \dt [e^{\mp i\frac{t}{\varepsilon}}(\Psi \pm i \Psi^\perp )]  = - v\cdot \nablah [e^{\mp i\frac{t}{\varepsilon}}(\Psi \pm i \Psi^\perp )] \\
    - w \dz [e^{\mp i\frac{t}{\varepsilon}}(\Psi \pm i \Psi^\perp )] 
    -  e^{\mp i\frac{t}{\varepsilon}}(N_2 \pm i N_2^\perp),\end{gathered} \\
    \label{eq:fast-wave-2}
    \text{and} \qquad
    \dt [e^{\mp i\frac{t}{\varepsilon}}(Z \pm i Z^\perp)] = - e^{\mp i\frac{t}{\varepsilon}}(N_3 \pm i N_3^\perp).
\end{gather}
From \eqref{eq:fast-wave-1} and \eqref{eq:fast-wave-2}, thanks to \eqref{est:uniform}, \eqref{uni-est:001}, and \eqref{uni-est:001-1}, it follows that 
\begin{equation}\label{uni-est:004}
    \begin{gathered}
    \sup_{{0\leq t \leq T}} \bigl(\norm{\dt [e^{\mp i\frac{t}{\varepsilon}}(\Psi(t) \pm i \Psi^\perp(t))]}{\Hnorm{1}} + \norm{\dt [e^{\mp i\frac{t}{\varepsilon}}(Z(t) \pm i Z^\perp(t))]}{\Hnorm{2}} \\
    + \norm{e^{\mp i\frac{t}{\varepsilon}}(\Psi(t) \pm i \Psi^\perp(t))}{\Hnorm{2}} + \norm{e^{\mp i\frac{t}{\varepsilon}}(Z(t) \pm i Z^\perp(t))}{\Hnorm{3}} \bigr) \\
    \leq C_\mrm{in}.\end{gathered}
\end{equation}
Therefore, by the Aubin compactness theorem \cite[Theorem 2.1]{Temam}, there exist
\begin{equation}\label{def:limit-variables-2}
    \begin{gathered}
    \psi_{p,\pm} \in L^\infty(0,T;H^2(\Omega)), \qquad\qquad z_{p,\pm} \in L^\infty(0,T;H^3(\Omega)),\\
    \dt \psi_{p,\pm} \in L^\infty(0,T;H^1(\Omega)), \qquad \text{and}\qquad \dt z_{p,\pm} \in L^\infty(0,T;H^2(\Omega)),
    \end{gathered}
\end{equation}
such that there exists a subsequence of $ \varepsilon $ that as $ \varepsilon \rightarrow 0^+ $,
\begin{align}
    \Psi_{\pm} 
    & \stackrel{\ast}{\rightharpoonup} \psi_{p,\pm} & \text{weak-$\ast$ in} & \qquad L^\infty(0,T;H^2(\Omega)), \label{cnvg:weak-3} \\
    Z_{\pm}
    & \stackrel{\ast}{\rightharpoonup} z_{p,\pm} & \text{weak-$\ast$ in} & \qquad L^\infty(0,T;H^3(\Omega)), \label{cnvg:weak-4} \\
    \dt \Psi_{\pm} 
    & \stackrel{\ast}{\rightharpoonup} \dt \psi_{p,\pm} & \text{weak-$\ast$ in} & \qquad L^\infty(0,T;H^1(\Omega)), \label{cnvg:weak-5} \\
    \dt Z_{\pm}
    & \stackrel{\ast}{\rightharpoonup} \dt z_{p,\pm} & \text{weak-$\ast$ in} & \qquad L^\infty(0,T;H^{2}(\Omega)), \label{cnvg:weak-5-2} \\
    \shortintertext{and}
    \Psi_{\pm}
    & \rightarrow \psi_{p,\pm} & \text{in} & \qquad C([0,T];H^1(\Omega)), \label{cnvg:strong-2}\\
    Z_{\pm} 
    & \rightarrow z_{p,\pm} & \text{in} & \qquad C([0,T];H^2(\Omega)), \label{cnvg:strong-3}
\end{align}
where
\begin{equation}\label{def:fast-variable}
    \Psi_{\pm} := e^{\mp i\frac{t}{\varepsilon}}(\Psi(t) \pm i \Psi^\perp(t)) \quad\text{and}\quad
    Z_{\pm} :=  e^{\mp i\frac{t}{\varepsilon}}(Z(t) \pm i Z^\perp(t)). 
\end{equation}
In particular, directly one can verify that
\begin{gather}
\begin{gathered}
    2 \Psi - (e^{i\frac{t}{\varepsilon}}\psi_{p,+} + e^{-i\frac{t}{\varepsilon}}\psi_{p,-}) = e^{i\frac{t}{\varepsilon}} (\Psi_+ - \psi_{p,+}) + e^{-i\frac{t}{\varepsilon}}(\Psi_- - \psi_{p,-}) \\
    \rightarrow 0 \qquad\text{in}\qquad L^\infty(0,T;H^1(\Omega)), \quad \text{as} ~ \varepsilon \rightarrow 0^+,
\end{gathered}\label{cnvg:strong-4}\\
\shortintertext{and}
\begin{gathered}
    2 Z - (e^{i\frac{t}{\varepsilon}} z_{p,+} + e^{-i\frac{t}{\varepsilon}}z_{p,-}) = e^{i\frac{t}{\varepsilon}} (Z_+ - z_{p,+}) + e^{-i\frac{t}{\varepsilon}}(Z_- - z_{p,-}) \\
    \rightarrow 0 \qquad\text{in}\qquad L^\infty(0,T;H^2(\Omega)), \quad \text{as} ~ \varepsilon \rightarrow 0^+.
\end{gathered}\label{cnvg:strong-5}
\end{gather}

To conclude this section, we write the fast-slow-error decomposition of $ v,w,\theta $. 
Let
\begin{equation}\label{def:fast-variables}
\begin{gathered}
    W_\pm :=  \dfrac{1}{2} \Delta_D^{-1}\dvh \psi_{p,\pm}, \qquad 
    \Theta_\pm :=  \dfrac{1}{2} \Delta_D^{-1}\crlh \psi_{p,\pm},
    \qquad \text{and}\\ 
    V_\pm :=  \dfrac{1}{2} \bigl(z_{p,\pm} - \nablah \Deltah^{-1}\dz \Delta_D^{-1} \dvh \psi_{p,\pm} - \nablah^\perp \Deltah^{-1}\dz \Delta_D^{-1} \crlh \psi_{p,\pm} \bigr). 
\end{gathered}
\end{equation}
Thanks to \eqref{def:limit-variables-2}, one has that
\begin{equation}\label{dmp:cmpt_w_tm_v} 
\begin{gathered}
    W_\pm, \Theta_\pm, V_\pm \in L^\infty(0,T;H^3(\Omega)) \\
    \text{and} \qquad \dt W_\pm, \dt \Theta_\pm, \dt V_\pm \in L^\infty(0,T;H^2(\Omega)).
\end{gathered}
\end{equation}
Repeating the exact calculation as in \eqref{eq:laplace_w_l}--\eqref{sol:v_l} leads to
\begin{align}
    \label{dmp:w-f-s-e}
    & \begin{aligned}
        w = & \Delta_D^{-1} \dvh \Psi =
        \underbrace{e^{i\frac{t}{\varepsilon}}W_+}_{=:w_{\mrm{fast},+}} + \underbrace{e^{-i\frac{t}{\varepsilon}}W_-}_{=: w_{\mrm{fast},-}}
        + w_\mrm{err},
    \end{aligned}\\
    \label{dmp:tm-f-s-e}
    & \begin{aligned}
        \theta = & \mathrm E_b(H_0,H_h) + \Delta_D^{-1}(\crlh \Psi + \dz \Phi - \Delta E_b(H_0,H_h)) \\
        = & \underbrace{\mathrm E_b(H_0,H_h) + \Delta_D^{-1}(\dz \Phi - \Delta E_b(H_0,H_h))}_{=:\theta_{\mrm{slow}}} 
        + \underbrace{e^{i\frac{t}{\varepsilon}}\Theta_+}_{=:\theta_{\mrm{fast},+}} 
        + \underbrace{e^{-i\frac{t}{\varepsilon}}\Theta_-}_{=:\theta_{\mrm{fast},-}}
        \\ & \qquad + \theta_\mrm{err},
    \end{aligned}\\
    \shortintertext{and}
    \label{dmp:v-f-s-e}
    & \begin{aligned}
        v = & Z - \nablah \Deltah^{-1} \dz (\Delta_D^{-1}\dvh \Psi) \\
        & + \nablah^\perp \Deltah^{-1} [\Phi - \dz \mathrm E_b(H_0,H_h) \\
        & \qquad\qquad\qquad - \dz \Delta_D^{-1}(\crlh \Psi + \dz \Phi - \Delta \mathrm E_b(H_0,H_h))] \\
        = & \underbrace{\substack{\nablah^\perp \Deltah^{-1} [\Phi - \dz \mathrm E_b(H_0,H_h) \\
        \qquad - \dz \Delta_D^{-1}(\dz \Phi - \Delta \mathrm E_b(H_0,H_h))]}}_{=:v_\mrm{slow}} 
        + \underbrace{e^{i\frac{t}{\varepsilon}}V_+}_{=:v_{\mrm{fast},+}} 
        + \underbrace{e^{-i\frac{t}{\varepsilon}}V_-}_{=:v_{\mrm{fast},-}}
        + v_\mrm{err},
    \end{aligned}
\end{align}
where, thanks to \eqref{uni-est:001}, \eqref{def:limit-variables-2}, \eqref{cnvg:strong-4}, and \eqref{cnvg:strong-5}, the error terms satisfy 
\begin{equation}\label{cnvg:strong-6}
    \begin{gathered}
        \sup_{0\leq t \leq T}\norm{v_\mrm{err}(t),w_\mrm{err}(t),\theta_\mrm{err}(t)}{\Hnorm{3}} \leq C_\mrm{in},\\
        v_\mrm{err}, ~ w_\mrm{err}, ~ \text{and} ~ \theta_\mrm{err} \rightarrow 0 \qquad \text{in}\qquad L^\infty(0,T;H^2(\Omega)), \quad \text{as} ~ \varepsilon \rightarrow 0^+.
    \end{gathered}
\end{equation}
In addition, thanks to \eqref{def:extension-to-interial},
\eqref{uni-est:003}, \eqref{cnvg:weak-1}, \eqref{cnvg:weak-1-2}, \eqref{cnvg:strong}, \eqref{cnvg:strong-1} and \eqref{def:slow-variables-2}, we have
\begin{equation}\label{cnvg:strong-7-1}
\begin{gathered}
\sup_{0\leq t\leq T}\norm{v_\mrm{slow},\theta_\mrm{slow}}{H^3(\Omega)} \leq C_\mrm{in} \qquad \text{and} \\
    \sup_{0\leq t \leq T} \norm{\dt v_\mrm{slow}(t), \dt \theta_\mrm{slow}(t)}{H^2(\Omega)} \leq C_\mrm{in}.
\end{gathered}
\end{equation}
Moreover, there exists a subsequence of $ \varepsilon $ that as $ \varepsilon \rightarrow 0^+ $, we also have
\begin{equation}\label{cnvg:strong-7}
        \begin{gathered}
        v_\mrm{slow}, \theta_\mrm{slow} \qquad \rightarrow \qquad v_p, \theta_p \qquad \text{in}  \qquad C(0,T;H^{2}(\Omega)), \\
        \text{and} \qquad v_\mrm{slow}, \theta_\mrm{slow} \qquad  \stackrel{\ast}{\rightharpoonup} \qquad v_p, \theta_p  \qquad  \text{weak-$\ast$ in} \qquad  L^\infty(0,T;H^3(\Omega)),\\
        \text{as} ~ \varepsilon \rightarrow 0^+.
        \end{gathered}
\end{equation}

\subsection{Convergence theory: Part 2, convergence of the nonlinearities }\label{sec:convergence_2}

In this section, we finish the convergence theory by investigating the convergence of the nonlinearities. 

\subsubsection*{Convergence of the slow waves \eqref{eq:pv}, \eqref{eq:bc_th_0}, and \eqref{eq:bc_th_h}}

First, we investigate $ N_1 $, defined in \eqref{def:nonln_pv}. Notice that $ N_1 $ is quadratic. substituting \eqref{dmp:w-f-s-e}--\eqref{dmp:v-f-s-e}, we write
\begin{align*}
     N_1 & =  \underbrace{\crlh v_\mrm{slow} \cdot \dvh v_\mrm{slow} + \dz v_\mrm{slow} \cdot \nablah \theta_\mrm{slow}}_{=:N_{1,\mrm{slow}}} \\
    & + \underbrace{
        \splitfrac
        {\crlh v_{\mrm{fast},\pm}\cdot \dvh v_{\mrm{fast},\mp} + \dz v_{\mrm{fast},\pm} \cdot \nablah^\perp w_{\mrm{fast},\mp}}
        { + \dz v_{\mrm{fast},\pm}\cdot\nablah \theta_{\mrm{fast},\mp} + \dz w_{\mrm{fast},\pm} \dz \theta_{\mrm{fast},\mp}}
    }_{=:N_{1,\mrm{res}}} \\
    & + \underbrace{
        \substack{
        \crlh v_{\mrm{slow}}\cdot \dvh v_{\mrm{fast},\pm} + \crlh v_{\mrm{fast},\pm}\cdot \dvh v_{\mrm{slow}} \\
        + \dz v_{\mrm{slow}} \cdot \nablah^\perp w_{\mrm{fast},\pm} + \dz v_{\mrm{slow}}\cdot\nablah \theta_{\mrm{fast},\pm}
        \\
         + \dz v_{\mrm{fast},\pm}\cdot\nablah \theta_{\mrm{slow}}+ \dz w_{\mrm{fast},\pm} \dz \theta_{\mrm{slow}} 
        }
    }_{=:N_{1,\mrm{fast},1}}\\
    & + \underbrace{
        \splitfrac
        {\crlh v_{\mrm{fast},\pm}\cdot \dvh v_{\mrm{fast},\pm} + \dz v_{\mrm{fast},\pm} \cdot \nablah^\perp w_{\mrm{fast},\pm}}
        { + \dz v_{\mrm{fast},\pm}\cdot\nablah \theta_{\mrm{fast},\pm} + \dz w_{\mrm{fast},\pm} \dz \theta_{\mrm{fast},\pm}}
    }_{=:N_{1,\mrm{fast},2}}\\
    & + \underbrace{\text{the rest terms}}_{=:N_{1,\mrm{err}}}.
\end{align*}
Then thanks to \eqref{est:uniform},
\eqref{dmp:cmpt_w_tm_v}, \eqref{cnvg:strong-6}, and \eqref{cnvg:strong-7}, we have, as $ \varepsilon \rightarrow 0^+ $,
\begin{align}
    N_{1,\mrm{slow}} & \rightarrow \crlh v_p \cdot \dvh v_p + \dz v_p \cdot \nablah \theta_p = 0 \quad \text{in} ~ C([0,T];H^1(\Omega)), \label{limit:N_1_slow} \\
    N_{1,\text{fast},1}, N_{1,\text{fast},2} & 
    \rightharpoonup 0 \qquad \text{weakly in} \quad L^p(0,T;H^{1}(\Omega)) \quad \forall p\in(1,\infty), \label{limit:N_1_fast_1} \\
    N_{1,\mrm{err}} & \rightarrow 0 \qquad\qquad\quad ~ \text{in} \quad L^\infty(0,T;H^1(\Omega)), \label{limit:N_1_fast_2}
\end{align}
and
\begin{equation}\label{limit:N_1_res}
    \begin{aligned}
        N_{1,\mrm{res}} & \rightarrow \crlh V_{\pm}\cdot \dvh V_{\mp} + \dz V_{\pm} \cdot \nablah^\perp W_{\mp}
        + \dz V_{\pm}\cdot\nablah \Theta_{\mp} + \dz W_{\pm} \dz \Theta_{\mp}\\
        \text{in} & \qquad L^\infty(0,T;H^2(\Omega)).
    \end{aligned}
\end{equation}

Consequently, as $ \varepsilon\rightarrow 0^+ $, in the sense of distribution, the limit of equation \eqref{eq:pv} is
\begin{equation}\label{limit:eq_pv}
    \begin{gathered}
        \dt \Phi_p + v_p \cdot \nablah \Phi_p + w_p \dz \Phi_p 
        + \crlh V_{\pm}\cdot \dvh V_{\mp} \\
        + \dz V_{\pm} \cdot \nablah^\perp W_{\mp}
        + \dz V_{\pm}\cdot\nablah \Theta_{\mp} + \dz W_{\pm} \dz \Theta_{\mp}
         = 0.
    \end{gathered}
\end{equation}
Here we have omitted the convergence of the advection terms, which is left to the reader. 

The limit equations of \eqref{eq:bc_th_0} and \eqref{eq:bc_th_h} follow similarly. The proof is left to the reader and we only state the result as follows:
\begin{gather}
    \dt H_{p,0} + v_p\vert_{z=0} \cdot \nablah H_{p,0} = 0, \label{limit:eq_bc_th_0}\\
    \dt H_{p,h} + v_p\vert_{z=h} \cdot \nablah H_{p,h} = 0. \label{limit:eq_bc_th_h}
\end{gather}

\bigskip

We remind the reader that $ v_p, w_p, \theta_p $ ($V_\pm, W_\pm, \Theta_\pm$, respectively) are determined by $ \Phi_p, H_{p,0}, H_{p,h} $ ($\psi_{p,\pm},z_{p,\pm})$, respectively), as in \eqref{def:slow-variables-2} (\eqref{def:fast-variables}, respectively). Therefore, the equations for $ \Phi_p $, $ H_{p,0}$, and $ H_{p,h} $, i.e., \eqref{limit:eq_pv}, \eqref{limit:eq_bc_th_0}, and \eqref{limit:eq_bc_th_h}, can be considered as the equations of $ v_p, w_p, \theta_p $, with source terms given by the resonances involving $ V_\pm $, $ W_\pm $, and $ \Theta_\pm $ (equivalently $ \psi_{p,\pm} $ and $ z_{p,\pm} $). To close the system, we will investigate the limit equations of \eqref{eq:fast-wave-1} and \eqref{eq:fast-wave-2} in the following.

\subsubsection*{Convergence of the fast waves \eqref{eq:fast-wave-1} and \eqref{eq:fast-wave-2}}

Using the notation of  \eqref{def:fast-variable}, \eqref{eq:fast-wave-1} and \eqref{eq:fast-wave-2} can be written as 
\begin{gather*}
    \dt \Psi_\pm + v\cdot \nablah \Psi_\pm + w \dz \Psi_\pm + e^{\mp i\frac{t}{\varepsilon}}(N_2 \pm i N_2^\perp) = 0 , \tag{\ref{eq:fast-wave-1}'} \\
    \dt Z_\pm + e^{\mp i \frac{t}{\varepsilon}}(N_3 \pm i N_3^\perp)=0. \tag{\ref{eq:fast-wave-2}'}
\end{gather*}
Thanks to \eqref{cnvg:weak-2} and \eqref{cnvg:weak-3}--\eqref{cnvg:strong-3}, we only need to investigate the limit of $ e^{\mp i \frac{t}{\varepsilon}}N_2 $ and $ e^{\mp i \frac{t}{\varepsilon}}N_3 $. 

Repeating the same arguments as for $ N_1 $, above, one can show that 
\begin{align*}
    e^{\mp i \frac{t}{\varepsilon}}N_{2} = & e^{\mp i\frac{t}{\varepsilon}} \bigl( (\nablah v_{\mrm{fast},\pm})^\top \nablah \theta_\mrm{slow} + (\nablah v_\mrm{slow})^\top \nablah \theta_{\mrm{fast},\pm} \bigr)^\perp \\
    & + e^{\mp i\frac{t}{\varepsilon}} \dz \theta_\mrm{slow}\cdot \nablah^\perp w_{\mrm{fast},\pm} + e^{\mp i\frac{t}{\varepsilon}} (\nablah v_\mrm{slow})^\top \nablah w_{\mrm{fast},\pm}\\
    & - e^{\mp i\frac{t}{\varepsilon}} (\dz v_{\mrm{fast},\pm}\cdot \nablah v_\mrm{slow} + \dz v_\mrm{slow} \cdot \nablah v_{\mrm{fast},\pm}) \\
    & - e^{\mp i\frac{t}{\varepsilon}} \dz w_{\mrm{fast},\pm} \dz v_{\mrm{slow}}
     + \underbrace{\text{the rest}}_{\mathclap{
     \rightharpoonup 0 \quad \text{in the sense of distribution}
     }}.
\end{align*}
After substituting \eqref{dmp:w-f-s-e}--\eqref{dmp:v-f-s-e} and sending $ \varepsilon \rightarrow 0^+ $, it follows that
\begin{equation}\label{def:nonln_psi}
    \begin{aligned}
        & e^{\mp i \frac{t}{\varepsilon}}N_{2} \rightharpoonup \bigl( (\nablah V_\pm)^\top \nablah \theta_p + (\nablah v_p)^\top \nablah \Theta_\pm  \bigr)^\perp \\ 
        & \qquad + \dz \theta_p \cdot \nablah^\perp W_\pm + (\nablah v_p)^\top \nablah W_\pm \\
        & \qquad - \dz V_\pm \cdot \nablah v_p - \dz v_p \cdot \nablah V_\pm - \dz W_\pm \dz v_p
        =: N_{\psi}\\
        & \qquad  \text{in}\quad L^p(0,T;H^1(\Omega))\quad \forall p \in (1,\infty).
    \end{aligned}
\end{equation}
Therefore, the limit of \eqref{eq:fast-wave-1} as $ \varepsilon \rightarrow 0^+ $ is
\begin{equation}\label{limit:eq_fast_wave_1}
    \dt \psi_{p,\pm} + v_p \cdot \nablah \psi_{p,\pm} + w_p \dz \psi_{p,\pm} + (N_\psi \pm i N_\psi^\perp) = 0.
\end{equation}

Last but not least, one has that
\begin{align*}
    e^{\mp i \frac{t}{\varepsilon}} N_3 = & e^{\mp i \frac{t}{\varepsilon}} \int_{\mathbb T^2} \dz (w_{\mrm{fast},\pm} v_\mrm{slow}) \,dxdy + \underbrace{\text{the rest}}_{\mathclap{\rightharpoonup 0 \quad \text{in the sense of distribution}}}, 
\end{align*}
and thus
\begin{equation}\label{def:nonln_Z}
    \begin{gathered}
        e^{\mp i \frac{t}{\varepsilon}} N_3 \rightharpoonup  \int_{\mathbb{T}^2} \dz (  W_\pm v_p) \,dxdy =: N_z \\
        \qquad \text{in}\quad L^p(0,T;H^1(\Omega))\quad \forall p \in (1,\infty).
    \end{gathered}
\end{equation}
Consequently, as $ \varepsilon \rightarrow 0^+ $, the limit of \eqref{eq:fast-wave-2} is 
\begin{equation}\label{limit:eq_fast_wave_2}
    \dt z_{p,\pm} + (N_z \pm i N_z^\perp) = 0.
\end{equation}

\subsubsection*{Conclusion}

The limit system for the slow limit variables $ \Phi_p, H_p $ and the fast limit variables $ \psi_{p,\pm}, z_{p,\pm} $ is then, from \eqref{limit:eq_pv}, \eqref{limit:eq_bc_th_0}, \eqref{limit:eq_bc_th_h}, \eqref{limit:eq_fast_wave_1}, and \eqref{limit:eq_fast_wave_2},
\begin{subequations}\label{system:limit}
    \begin{gather}
        \dt \Phi_p + v_p \cdot \nablah \Phi_p  +  N_\Phi  = 0,
        \label{limit:phi}\\
        \dt H_{p,0} + v_p\vert_{z=0} \cdot \nablah H_{p,0} = 0, \label{limit:bc_1}\\
        \dt H_{p,h} + v_p\vert_{z=h} \cdot \nablah H_{p,h} = 0, \label{limit:bc-2} \\
        \dt \psi_{p,\pm} + v_p \cdot \nablah \psi_{p,\pm} + (N_\psi \pm i N_\psi^\perp) = 0, \label{limit:psi}\\
        \dt z_{p,\pm} + (N_z \pm i N_z^\perp) = 0, \label{limit:z}
    \end{gather}
    where
    \begin{gather}
        \label{limit:pv_nonlinearity}
        N_\Phi := \crlh V_{\pm}\cdot \dvh V_{\mp}
        + \dz V_{\pm} \cdot \nablah^\perp W_{\mp}
        + \dz V_{\pm}\cdot\nablah \Theta_{\mp} + \dz W_{\pm} \dz \Theta_{\mp},
    \end{gather}
    and $ N_\psi $ and $N_z$  are defined in \eqref{def:nonln_psi} and \eqref{def:nonln_Z}, above, respectively. This finishes the proof of Theorem \ref{thm:convergence}.
\end{subequations}

\section*{Acknowledgements}

The authors would like to thank the Isaac Newton Institute for Mathematical Sciences, Cambridge, for support and hospitality during the programme ``Mathematical aspects of turbulence: where do we stand?'' (2022), where part of the work on this paper was undertaken. This work was supported in part by EPSRC grant no EP/R014604/1. X.L.'s work was partially supported by a grant from the Simons Foundation, during his visit to the Isaac Newton Institute for Mathematical Sciences. The research of EST has benefited from the inspiring environment of the CRC 1114 ``Scaling Cascades in Complex Systems'', Project Number
235221301, Project C06, funded by Deutsche Forschungsgemeinschaft (DFG).

\bibliographystyle{siam}

\end{document}